\documentclass[12pt]{amsart}
\usepackage{amssymb}

\let\frak\mathfrak
\let\Bbb\mathbb

\def\>{\relax\ifmmode\mskip.666667\thinmuskip\relax\else\kern.111111em\fi}
\def\<{\relax\ifmmode\mskip-.333333\thinmuskip\relax\else\kern-.0555556em\fi}
\def\vsk#1>{\vskip#1\baselineskip}
\def\vv#1>{\vadjust{\vsk#1>}\ignorespaces}
\def\vvn#1>{\vadjust{\nobreak\vsk#1>\nobreak}\ignorespaces}

\let\alb\allowbreak

\def\sskip{\par\vsk.2>}
\let\Medskip\medskip
\def\medskip{\par\Medskip}
\let\Bigskip\bigskip
\def\bigskip{\par\Bigskip}

\let\Maketitle\maketitle
\def\maketitle{\Maketitle\thispagestyle{empty}\let\maketitle\empty}

\newtheorem{thm}{Theorem}[section]
\newtheorem{cor}[thm]{Corollary}
\newtheorem{lem}[thm]{Lemma}

\numberwithin{equation}{section}

\theoremstyle{definition}
\newtheorem*{rem}{Remark}
\newtheorem*{example}{Example}

\let\mc\mathcal
\let\nc\newcommand

\nc{\on}{\operatorname}
\nc{\Z}{{\mathbb Z}}
\nc{\C}{{\mathbb C}}
\nc{\N}{{\mathbb N}}
\nc{\pone}{{\mathbb C}{\mathbb P}^1}
\nc{\arr}{\rightarrow}
\nc{\larr}{\longrightarrow}
\nc{\al}{\alpha}
\nc{\W}{{\mc W}}
\nc{\la}{\lambda}
\nc{\su}{\widehat{{\mathfrak sl}}_2}
\nc{\g}{{\mathfrak g}}
\nc{\h}{{\mathfrak h}}
\nc{\m}{{\mathfrak m}}
\nc{\n}{{\mathfrak n}}
\nc{\Gm}{\Gamma}
\nc{\La}{\Lambda}
\nc{\gl}{\widehat{\mathfrak{gl}_2}}
\nc{\bi}{\bibitem}
\nc{\om}{\omega}
\nc{\Res}{\on{Res}}
\nc{\gm}{\gamma}
\nc{\Om}{\Omega}

\def\fratop{\genfrac{}{}{0pt}1}

\def\ev{\mbox{\sl ev}}

\def\End{\on{End}}

\def\Gr{\on{Gr}}
\def\Res{\on{Res}}
\def\rdet{\on{rdet}}
\def\Wr{\on{Wr}}

\def\Llb{{\bs\La,\bs\la\>,\bs b}}

\def\B{{\mc B}}
\def\D{{\mc D}}

\def\M{{\mc M}}
\def\O{{\mc O}}

\def\flati{\def\=##1{\rlap##1\hphantom b)}}

\let\dl\delta
\let\Dl\Delta

\let\si\sigma

\let\der\partial

\let\ge\geqslant

\let\le\leqslant
\let\leq\leqslant

\nc{\gln}{\mathfrak{gl}_N}
\nc{\sln}{\mathfrak{sl}_N}

\def\glnt{\glN[t]}
\def\Ugln{U(\glN)}
\def\Uglnt{U(\glNt)}

\nc{\glN}{\mathfrak{gl}_N}
\nc{\slN}{\mathfrak{sl}_N}

\def\glNt{\glN[t]}

\def\beq{\begin{equation}}
\def\eeq{\end{equation}}
\def\be{\begin{equation*}}
\def\ee{\end{equation*}}

\nc{\bean}{\begin{eqnarray}}
\nc{\eean}{\end{eqnarray}}
\nc{\bea}{\begin{eqnarray*}}
\nc{\eea}{\end{eqnarray*}}
\nc{\bs}{\boldsymbol}
\nc{\Ref}[1]{{\rm(\ref{#1})}}

\nc{\Wl}{\Wr_{\bs\la}}
\nc{\Ol}{\Om_{\bs\la}}
\nc{\Mla}{\M_{\bs\la,\bs a}}
\nc{\p}{\partial}
\nc{\Bl}{\B_{\bs\la}}
\nc{\Bla}{\B_{\bs\la,\bs a}}
\nc{\Ola}{\O_{\bs\la,\bs a}}

\nc{\OL}{\mc O_{\bs\La,\bs\la,\bs b}}

\nc{\Bm}{\mc B_{m^{\bs\La}_{\bs\la}}}
\nc{\R}{\Bbb R}

\nc{\s}{{\rm sing}}
\nc{\Oll}{{\Omega_{\bs\la}}}

\def\sing{{\rm\>sing}}

\def\fratop{\genfrac{}{}{0pt}1}
\def\satop#1#2{\fratop{\scriptstyle#1}{\scriptstyle#2}}

\nc{\Tee}{\mathcal S}
\nc{\un}{U(\n_-)}

\nc{\ep}{\epsilon}

\begin{document}

\hrule width0pt
\vsk->

\title[On separation of variables for {\small\;$\glN$} Gaudin model]
{On separation of variables and completeness of\\[3pt] the Bethe ansatz
for quantum {\large\,$\glN$} Gaudin model}

\author[E.\,Mukhin, V.\,Tarasov, and \>A.\,Varchenko]
{E.\,Mukhin$\>^*$, V.\,Tarasov$\>^\star$, and \>A.\,Varchenko$\>^\diamond$}

\maketitle

\begin{center}
{\it $\kern-.4em^{*,\star}\<$Department of Mathematical Sciences,
Indiana University\,--\>Purdue University Indianapolis\kern-.4em\\
402 North Blackford St, Indianapolis, IN 46202-3216, USA\/}

\medskip
{\it $^\star\<$St.\,Petersburg Branch of Steklov Mathematical Institute\\
Fontanka 27, St.\,Petersburg, 191023, Russia\/}

\medskip
{\it $^\diamond\<$Department of Mathematics, University of North Carolina
at Chapel Hill\\ Chapel Hill, NC 27599-3250, USA\/}
\end{center}

{\let\thefootnote\relax
\footnotetext{\vsk-.8>\noindent
$^*$\,Supported in part by NSF grant DMS-0601005\\
$^\star$\,Supported in part by RFFI grant 05-01-00922\\
$^\diamond$\,Supported in part by NSF grant DMS-0555327}}

\medskip
\begin{abstract}
In this note, we discuss implications of the results obtained in \cite{MTV4}.
It was shown there that eigenvectors of the Bethe algebra of the quantum $\glN$
Gaudin model are in a one-to-one correspondence with Fuchsian differential
operators with polynomial kernel. Here, we interpret this fact as a separation
of variables in the $\glN$ Gaudin model. Having a Fuchsian differential
operator with polynomial kernel, we construct the corresponding eigenvector of
the Bethe algebra. It was shown in \cite{MTV4} that the Bethe algebra has
simple spectrum if the evaluation parameters of the Gaudin model are generic.
In that case, our Bethe ansatz construction produces an eigenbasis of the Bethe
algebra.
\end{abstract}

\setcounter{footnote}{0}
\renewcommand{\thefootnote}{\arabic{footnote}}

\section{Introduction}
Generally speaking, separation of variables in a quantum integrable
model is a reduction of a multidimensional spectral problem to a
suitable one-dimensional problem. For example, the famous Sklyanin's
separation of variables for the $\frak{gl}_2$ Gaudin model \cite{Sk} is
a reduction of the diagonalization problem of the Gaudin Hamiltonians,
acting on a tensor product of $\frak{gl}_2$-modules, to the problem of
finding a second order Fuchsian differential operator with polynomial
kernel and prescribed singularities. Having such a differential operator,
Sklyanin constructs an eigenvector of the Hamiltonians.

It has been proved recently in \cite{MTV4} that the eigenvectors of the Bethe
algebra of the $\glN$ Gaudin model are in a bijective correspondence with
$N$-th order Fuchsian differential operators with polynomial kernel and
prescribed singularities. This reduces the multidimensional problem of
the diagonalization of the Bethe algebra to the one-dimensional problem of
finding the corresponding Fuchsian differential operators. In that respect,
``the variables are separated''.

Having an eigenvector of the Bethe algebra, one has an effective way to
construct the corresponding Fuchsian operator, see~\cite{MTV2}, \cite{MTV4},
and Theorem~\ref{thm eigen to Fuchs}. In the opposite direction, the assignment
of an eigenvector to a Fuchsian operator is not explicit in \cite{MTV4}.
In this note, having a Fuchsian differential operator with polynomial kernel,
we construct the corresponding eigenvector of the Bethe algebra.
Our construction of an eigenvector from a differential operator can be viewed
as a (generalized) Bethe ansatz construction, cf.~\cite{Ba}, \cite{RV},
\cite{MV1}, \cite{MV2}.

It has been proved in \cite{MTV4}, that the action of the Bethe algebra on
a tensor product of irreducible finite-dimensional evaluation $\glNt$-modules
has simple spectrum provided the evaluation points are generic. In that case,
our construction of eigenvectors of the Bethe algebra produces an eigenbasis
of the Bethe algebra, thus showing the completeness of the Bethe ansatz.

\medskip
The authors thank J.\,Harnad for useful discussions.

\section{Eigenvectors of Bethe algebra}
\label{alg sec}
\subsection{Lie algebra $\gln$}
Let $e_{ij}$, $i,j=1,\dots,N$, be the standard generators of the Lie algebra
$\gln$ satisfying the relations
$[e_{ij},e_{sk}]=\dl_{js}e_{ik}-\dl_{ik}e_{sj}$.

Let $M$ be a $\gln$-module. A vector $v\in M$ has weight
$\bs\la=(\la_1,\dots,\la_N)\in\C^N$ if $e_{ii}v=\la_iv$ for $i=1,\dots,N$.
A vector $v$ is called {\it singular\/} if $e_{ij}v=0$ for $1\le i<j\le N$.

We denote by $(M)_{\bs\la}$ the subspace of $M$ of weight $\bs\la$,
by $(M)^\sing$ the subspace of $M$ of all singular vectors and by
$(M)_{\bs\la}^\sing$ the subspace of $M$ of all singular vectors
of weight $\bs\la$.

Denote by $L_{\bs\la}$ the irreducible finite-dimensional $\gln$-module with
highest weight $\bs\la$. Any finite-dimensional $\gln$-module $M$ is isomorphic
to the direct sum $\bigoplus_{\bs\la}L_{\bs\la}\otimes (M)_{\bs\la}^\sing$,
where the spaces $(M)_{\bs\la}^\sing$ are considered as trivial $\gln$-modules.

The $\gln$-module $L_{(1,0,\dots,0)}$ is the standard $N$-dimensional vector
representation of $\gln$. We denote it by $V$. We choose a highest weight
vector in $V$ and denote it by $v_+$.

A $\gln$-module $M$ is called polynomial if it is isomorphic to a submodule of
$V^{\otimes n}$ for some $n$.

A sequence of integers $\bs\la=(\la_1,\dots,\la_N)$ such that
$\la_1\ge\la_2\ge\dots\ge\la_N\ge0$ is called a {\it partition with at most
$N$ parts\/}. Set $|\bs\la|=\sum_{i=1}^N\la_i$. Then it is said that $\bs\la$
is a partition of $|\bs\la|$.

\sskip
The $\gln$-module $V^{\otimes n}$ contains the module $L_{\bs\la}$
if and only if $\bs\la$ is a partition of $n$ with at most $N$ parts.

\subsection{Current algebra $\glnt$}
Let $\glnt=\gln\otimes\C[t]$ be the Lie algebra of $\gln$-valued polynomials
with the pointwise commutator. We call it the {\it current algebra\/}.
We identify the Lie algebra $\gln$ with the subalgebra $\gln\otimes1$
of constant polynomials in $\glnt$. Hence, any $\glnt$-module has the canonical
structure of a $\gln$-module.

It is convenient to collect elements of $\glnt$ in generating series
of a formal variable $u$. For $g\in\gln$, set
\vvn-.1>
\be
g(u)=\sum_{s=0}^\infty (g\otimes t^s)u^{-s-1}.
\vv.4>
\ee

For each $a\in\C$, there exists an automorphism $\rho_a$ of $\glnt$,
\;$\rho_a:g(u)\mapsto g(u-a)$. Given a $\glnt$-module $M$, we denote by $M(a)$
the pull-back of $M$ through the automorphism $\rho_a$. As $\gln$-modules,
$M$ and $M(a)$ are isomorphic by the identity map.

We have the evaluation homomorphism,
${\ev:\glnt\to\gln}$, \;${\ev:g(u) \mapsto g\>u^{-1}}$.
Its restriction to the subalgebra $\gln\subset\glnt$ is the identity map.
For any $\gln$-module $M$, we denote by the same letter the $\glnt$-module,
obtained by pulling $M$ back through the evaluation homomorphism. For each
$a\in\C$, the $\glnt$-module $M(a)$ is called an {\it evaluation module\/}.

\subsection{Bethe algebra}
\label{secbethe}
Given an $N\times N$ matrix $A$ with possibly noncommuting entries $a_{ij}$,
we define its {\it row determinant\/} to be
\vvn.3>
\bea
\rdet A\,=
\sum_{\;\si\in S_N\!} (-1)^\si\,a_{1\si(1)}a_{2\si(2)}\dots a_{N\si(N)}\,.
\vv.2>
\eea

Let $\der$ be the operator of differentiation in variable $u$.
Define the {\it universal differential operator\/} $\D^\B$ by
\vvn-.6>
\be
\D^\B=\,\rdet\left( \begin{matrix}
\der-e_{11}(u) & -\>e_{21}(u)& \dots & -\>e_{N1}(u)\\
-\>e_{12}(u) &\der-e_{22}(u)& \dots & -\>e_{N2}(u)\\
\dots & \dots &\dots &\dots \\
-\>e_{1N}(u) & -\>e_{2N}(u)& \dots & \der-e_{NN}(u)
\end{matrix}\right).
\vv.2>
\ee
It is a differential operator in variable $u$, whose coefficients are
formal power series in $u^{-1}$ with coefficients in $\Uglnt$,
\vvn-.3>
\bea
\D^\B=\,\der^N+\sum_{i=1}^N\,B_i(u)\,\der^{N-i}\>,
\vv-.5>
\eea
where
\vvn-.3>
\bea
B_i(u)\,=\,\sum_{j=i}^\infty B_{ij}\>u^{-j}\,,
\eea
and $B_{ij}\in\Uglnt$, \,$i=1,\dots,N$, \,$j\in\Z_{\ge i}\>$.
We call the unital subalgebra of $\Uglnt$ generated by $B_{ij}$,
\,$i=1,\dots,N$, \,$j\in\Z_{\ge i}\>$, the {\it Bethe algebra\/}
and denote it by $\B$.

\sskip
By \cite{T}, \cite{MTV1}, the algebra $\B$ is commutative,
and the algebra $\B$ commutes with the subalgebra $\Ugln\subset \Uglnt$.

\sskip
As a subalgebra of $\Uglnt$, the algebra $\B$ acts on any $\glnt$-module $M$.
Since $\B$ commutes with $\Ugln$, it preserves the subspace of singular vectors
$(M)^\sing$ as well as weight subspaces of $M$. Therefore, the subspace
$(M)_{\bs\la}^\sing$ is $\B$-invariant for any weight $\bs\la$.

\sskip
Let $\bs\la^{(1)},\dots,\bs\la^{(k)},\>\bs\la$ \,be partitions with at most $N$
parts, and $b_1,\dots,b_k$ distinct complex numbers. We are interested
in the action of the Bethe algebra $\B$ on the tensor product of evaluation
modules $\otimes_{s=1}^kL_{\bs\la^{(s)}}(b_s)$, more precisely, on the subspace
$(\otimes_{s=1}^kL_{\bs\la^{(s)}}(b_s))_{\bs\la}^\sing$.

\sskip
Note that the subspace
$(\otimes_{s=1}^kL_{\bs\la^{(s)}}(b_s))_{\bs\la}^\sing$ is zero-dimensional
unless $|\bs\la|=\sum_{s=1}^k|\bs\la^{(s)}|$.

\subsection{Fuchsian differential operators and eigenvectors of Bethe algebra}
\label{sec Fuch oper}
Denote $\bs\La=(\bs\la^{(1)},\dots,\bs\la^{(k)})$ and $\bs b=(b_1,\dots,b_k)$.
Let $\Dl_\Llb$ be the set of all monic Fuchsian differential operators of
order $N$,
\bea
\D\,=\,\der^N+\>\sum_{i=1}^N\,h_i^\D(u)\,\der^{N-i}\,,
\eea
with the following properties.

\begin{enumerate}
\flati
\item[\=a]
The singular points of $\D$ are at\/ $b_1,\dots,b_k$ and $\infty$ only.
\item[\=b]
The exponents of $\D$ at\/ \>$b_s$\>, \,$s=1,\dots,k$, are
equal to$\,\la_N^{(s)},\,\la_{N-1}^{(s)}+1,\alb\,\dots\,,\la_1^{(s)}+N-1\,$.
\item[\=c]
The exponents of $\D$ at $\infty$ are equal to
$\>1-N-\la_1,\,2-N-\la_2,\,\dots\,,-\>\la_N\,$.
\item[\=d]
The kernel of the operator $\D$ consists of polynomials only.
\end{enumerate}

\noindent
Note that the set $\Dl_\Llb$ is empty
unless $|\bs\la|=\sum_{s=1}^k|\bs\la^{(s)}|$.

\medskip
Let $M$ be a $\glnt$-module and $v$ an eigenvector of the Bethe algebra
$\B\subset\Uglnt$ acting on $M$. Then for any coefficient $B_i(u)$ of
the universal differential operator $\D^\B$ we have $B_i(u)v=h_i(u)v$,
where $h_i(u)$ is a scalar series. We call the scalar differential operator
\bea
\D^\B_v\>=\,\der^N+\>\sum_{i=1}^N\>h_i(u)\,\der^{N-i}
\eea
the {\it differential operator associated with the eigenvector\/} $v$.

\begin{thm}
\label{thm eigen to Fuchs}
Let $v \in (\otimes_{s=1}^kL_{\bs\la^{(s)}}(b_s))_{\bs\la}^\sing$
be an eigenvector of the Bethe algebra, then
$\D^\B_v \in \Dl_\Llb$. Moreover, the assignment $v \mapsto \D^\B_v$ is
a bijective correspondence between the set of
eigenvectors of the action of the Bethe algebra on
$(\otimes_{s=1}^kL_{\bs\la^{(s)}}(b_s))_{\bs\la}^\sing$
(considered up to multiplication by nonzero numbers) and the set $\Dl_\Llb$.
\end{thm}

The first statement is Theorem 4.1 in \cite{MTV2}, cf.~\cite{MTV3}.
The second statement is Theorem 7.1 in \cite{MTV4}.

\medskip
The goal of this note is to construct the inverse bijection.

\section{Schubert cell and universal weight function}

\subsection{The cell $\Om_{\bs\la}$}
\label{Ominfty}
Let $N,d\in\Z_{>0}$, $N\leq d$. Let $\C_d[u]$ be the space of
polynomials in $u$ of degree less than $d$. We have $\dim \C_d[u]=d$.
Let $\Gr(N,d)$ be the Grassmannian of all $N$-dimensional subspaces in
$\C_d[u]$.

\medskip
Given a partition $\bs\la=(\la_1,\dots,\la_N)$ such that $\la_1\leq d-N$,
introduce a sequence
\vvn.2>
\be
P\,=\,\{d_1>d_2>\dots>d_N\}\,,\qquad d_i=\la_i+N-i\,,
\ee
and
denote by $\Omega_{\bs\la}$ the subset of $ \Gr(N,d)$ consisting of all
$N$-dimensional subspaces
$X\subset\C_d[u]$ such that
for every $i=1,\dots,N$, the subspace $X$ contains a polynomial
of degree $d_i$.

In other words, $\Om_{\bs\la}$ consists
of subspaces $X\subset\C_d[u]$ with a basis
$\{f_1(u),\dots,f_N(u)\}$ of the form
\vvn-.5>
\bea
\label{basis}
f_i(u)=u^{d_i}+\sum_{j=1,\ d_i-j\not\in P}^{d_i}f_{ij}u^{d_i-j}.
\vv.2>
\eea
For a given $X\in\Om_{\bs\la}$, such a basis is unique. The basis
$\{f_1(u),\dots,f_N(u)\}$ will be called the flag basis of the subspace $X$.

The set $\Om_{\bs\la}$ is a (Schubert) cell isomorphic to an affine space
of dimension $|\bs\la|$ with coordinate functions $f_{ij}$.

\sskip
For $X\in\Oll$, we denote by $\D_X$ the monic scalar differential operator of
order $N$ with kernel $X$.
We call $\D_X$ the {\it differential operator associated with\/} $X$.

\medskip

\subsection{Generic points of $\Om_{\bs\la}$}
\label{sec generic pts}
For $g_1,\dots,g_l \in \C[u]$, introduce the Wronskian by the formula
\vvn-.4>
\be
\Wr(g_1(u),\dots,g_l(u))\,=\,
\det\left(\begin{matrix} g_1(u) & g_1'(u) &\dots & g_1^{(l-1)}(u) \\
g_2(u) & g_2'(u) &\dots & g_2^{(l-1)}(u) \\ \dots & \dots &\dots & \dots \\
g_l(u) & g_l'(u) &\dots & g_l^{(l-1)}(u)
\end{matrix}\right).
\ee

For $X\in\Oll$\>, let $\{f_1(u),\dots,f_N(u)\}$ be the flag basis of $X$.
Introduce the polynomials
$\{y_0(u)\>,\,y_{1}(u)\>,\,\dots\,,\,y_{N-1}(u)\}$, \,by the formula
\vvn.3>
\be
y_a(u)\!\!\prod_{a<i<j\leq N} (\la_i-\la_j)\,
\,=\,\Wr(f_{a+1}(u),\dots,f_{N}(u))\,,\qquad a=0,\dots,N\,.
\ee
Set
\vvn-.6>
\beq
\label{sequence l}
l_a\,=\sum_{b=a+1}^N \la_b\,,\qquad a=0,\dots,N\,.
\eeq
Clearly, $l_0=|\bs\la|$ \,and \,$l_N=0$.

\sskip
For each $a=0,\dots,N-1$, the polynomial $y_a(u)$ is a monic polynomial of
degree $l_a$.
Denote $t_1^{(a)},\dots,t_{l_a}^{(a)}$ the roots of the polynomial $y_a(u)$,
and
\beq
\label{t z}
\bs t_X\,=\,(t_1^{(0)},\dots,t_{l_0}^{(0)},\,\dots\,,
t_1^{(N-1)},\dots,t_{l_{N-1}}^{(N-1)})\,.
\eeq
We say that $\bs t_X$ are the {\it root coordinates} of $X$.

We say that $X\in\Oll$ is {\it generic} if all roots of the polynomials
$y_0(u)\>,\,y_1(u)\>,\,\,\dots\,,\,y_{N-1}(u)$ are simple and for each
$a=1,\dots,N-1$, the polynomials $y_{a-1}(u)$ and $y_a(u)$ do not have common
roots.

\medskip

If $X$ is generic, then the root coordinates $\bs t_X$ satisfy the Bethe ansatz
equations \cite{MV1}:
\bea
\label{BAE Gaudin}
\sum_{j'=1}^{l_{a-1}}\frac 1{t^{(a)}_j - t^{(a-1)}_{j'}}\;-\,
\sum_{\satop{j'=1}{j'\neq j}}^{l_a}\frac 2{t^{(a)}_j - t^{(a)}_{j'}}\;+\,
\sum_{j'=1}^{l_{a+1}}\frac 1{t^{(a)}_j - t^{(a+1)}_{j'}}\;=\,0\,.
\vv-.2>
\eea
Here the equations are labeled by $a=1,\dots,N-1$, \,$j=1,\dots,l_a$.

\sskip
Conversely, if \,$\bs t\,=\,(t_1^{(0)},\dots,t_{l_0}^{(0)},\,\dots\,,
t_1^{(N-1)},\dots,t_{l_{N-1}}^{(N-1)})$ satisfy the Bethe ansatz equations,
then there exists a unique $X\in \Oll$ such that $X$ is generic and $\bs t$ are
its root coordinates, see \Ref{t z}. This $X$ is determined by the following
construction, see~\cite{MV1}. Set
\be
\chi^a(u, \bs t)\,=\,\sum_{j=1}^{l_{a-1}}\,\frac1{u-t^{(a-1)}_j}\;-\,
\sum_{i=1}^{l_a}\,\frac 1 {u- t^{(a)}_j}\;,\qquad a=1,\dots,N\,.
\ee
Then
\vvn-.5>
\be
\D_X\,=\,\bigl(\der -\chi^1(u,\bs t)\bigr)\,\dots\,
\bigl(\der-\chi^N(u,\bs t)\bigr)\,.
\ee

\begin{lem}
\label{lem on generic pts}
Generic points form a Zariski open subset of\/ $\Oll$.
\end{lem}

The lemma follows, for example, from part (i) of Theorem 6.1 in \cite{MV2}.

\subsection{Universal weight function}
Let $\bs\la$ be a partition with at most $N$ parts. Let $l_0,\dots,l_N$
be the numbers defined in \Ref{sequence l}. Denote $n=l_0$\>,
\,$l=l_1+\dots+l_{N-1}$ \>and \,$\bs l=(l_1,\dots,l_{N-1})$.

\sskip
Consider the weight subspace $(V^{\otimes n})_{\bs\la}$ of the $n$-th tensor
power of the vector representation of $\glN$ and the space $\C^{l+n}$ with
coordinates \,$\bs t\,=\,(t_1^{(0)},\dots,t_{l_0}^{(0)},\,\dots\,,
t_1^{(N-1)},\dots,t_{l_{N-1}}^{(N-1)})$.

In this section we remind the construction of a rational map
$\omega:\C^{l+n}\to (V^{\otimes n})_{\bs\la}$,
called the {\it universal weight function}, see~\cite{SV}.

A basis of $V^{\otimes n}$ is formed by the vectors
\vvn.1>
\be
e_J\>v\,=\,e_{j_1,1}\>v_+\otimes \dots\otimes e_{j_n,1}\>v_+\,,
\vv.2>
\ee
where $J=(j_1,\dots,j_n)$ and $1\leq j_s\leq N$ for $s=1,\dots,N$. A basis
of $(V^{\otimes n})_{\bs\la}$ is formed by the vectors $e_J\>v$ such that
$\#\{s\ |\ j_s>i\}\,=\,l_i$ for every $i=1,\dots,N-1$.
Such a $J$ will be called $\bs l$-admissible.

The universal weight function has the form
\vvn.3>
\be
\omega(\bs t)\,=\,\sum_J\,\omega_J(\bs t)\,e_Jv
\ee
where the sum is over the set of all $\bs l$-admissible $J$,
and the function $\omega_J(\bs t)$ is defined below.

\sskip
For an admissible $J$, \>define \,$S(J)=\{s\ |\ j_s>1 \}$\>,
and for \>$i=1,\ldots,N-1$, define
\vvn.3>
\be
S_i(J)\,=\,\{\,s\ |\ 1\le s\le n\,,\ \ 1\le i<j_s\,\}\,.
\vv-.2>
\ee
Then $|\>S_i(J)\>|\,=\,l_i$.

\sskip
Let $B(J)$ be the set of sequences \,$\bs\beta=(\beta_1,\dots,\beta_{N-1})$
of bijections \,$\beta_i:S_i(J)\to\{1,\dots,l_i\}$, \>$i=1,\dots,N-1$.
Then $|\>B(J)\>|\>=\prod_{a=1}^{N-1}l_a!$~.

\sskip
For $s\in S(J)$ and $\bs\beta\in B(J)$, introduce the rational function
\be
\omega_{s,\bs\beta}(\bs t)\,=\,\frac1{t^{(1)}_{\beta_1(s)}-t^{(0)}_s}\;
\prod_{i=2}^{j_1-1}\frac1{t^{(i)}_{\beta_i(s)}-t^{(i-1)}_{\beta_{i-1}(s)}}\
\ee
and define
\be
\omega_J(\bs t)\,=\,
\sum_{\bs\beta\in B(J)}\,\prod_{s\in S(J)}\,\omega_{s,\bs\beta} \ .
\ee

\begin{example}
Let $n=2$ \,and \,$\bs l=(1,1,0,\dots,0)$. Then
\vvn.5>
\be
\omega(\bs t)\,=\,\frac 1{(t_1^{(2)}-t_1^{(1)})\>
(t_1^{(1)}-t_1^{(0)})}\ e_{3,1}v_+\otimes v_+\,+\;
\frac 1{(t_1^{(2)}-t_1^{(1)})\>(t_1^{(1)}-t_2^{(0)})}\ v_+\otimes e_{3,1}v_+\;.
\vv.5>
\ee
\end{example}

\begin{thm}
\label{thm X to Vn}
Let $X\in \Oll$ be a generic point with root coordinates $\bs t_X$.
Consider the value $\omega(\bs t_X)$ of the universal weight function
$\omega:\C^{l+n}\to(V^{\otimes n})_{\bs\la}$ at\/ $\bs t_X$. Consider
$V^{\otimes n}$ as the $\glNt$-module $\otimes_{s=1}^nV(t_s^{(0)})$\>.
Then
\begin{enumerate}
\item[(i)]
The vector $\omega(\bs t_X)$ belongs to $(V^{\otimes n})_{\bs\la}^\sing$.

\item[(ii)] The vector $\omega(\bs t_X)$ is an eigenvector of the Bethe algebra
$\B$, acting on $\otimes_{s=1}^nV(t_s^{(0)})$. Moreover,
$\D^{\mc B}_{\omega(\bs t_X)}=\D_X$, where $\D^{\mc B}_{\omega(\bs t_X)}$
and $\D_X$ are the differential operators associated with the eigenvector
$\omega(\bs t_X)$ and the point $X\in\Oll$, respectively.
\end{enumerate}
\end{thm}

Part (i) is proved in~\cite{Ba} and~\cite{RV}. Part (i) also follows directly
from Theorem~6.16.2 in~\cite{SV}. Part (ii) is proved in \cite{MTV1}.

\begin{rem}
For a generic point $X\in \Oll$ the differential operator
$\D_X$ has the following properties.
\begin{enumerate}
\flati
\item[\=e]
The singular points of $\D_X$ are at\/ $t_1^{(0)},\dots,t_n^{(0)}$ and $\infty$
only;
\item[\=f]
The exponents of $\D_X$ at\/ \>$t_s^{(0)}$\>, \,$s=1,\dots,n$, are equal
to $\,0\>,\,1,\,\dots\,,N-2, N\,$;
\item[\=g]
The exponents of $\D_X$ at $\infty$ are equal to
$\>1-N-\la_1,\,2-N-\la_2,\,\dots\,,-\>\la_N\,$;
\item[\=h]
The kernel of the operator $\D_X$ consists of polynomials only.
\end{enumerate}

On the other hand, Theorem \ref{thm eigen to Fuchs}, applied to
the $\glNt$-module $\otimes_{s=1}^nV(t_s^{(0)})$, yields that
for any eigenvector $v$ of the Bethe algebra $\B$, acting on
$(\otimes_{s=1}^nV(t_s^{(0)}))_{\bs\la}^\sing$,
the differential operator $\D^\B_v$ has properties \,e)\,--\,h).

\sskip
Therefore, the universal weight function and the assignment
$\D_X\mapsto X\mapsto\omega(\bs t_X)$ allows us to reverse the correspondence
$v\mapsto\D^\B_v$ of Theorem~\ref{thm eigen to Fuchs} for the case of
the $\glNt$-module $\otimes_{s=1}^nV(t_s^{(0)})$ under the condition that
$X\in\Oll$ is generic. Our goal is to generalize this construction to the case
of a $\glNt$-module $\otimes_{s=1}^kL_{\bs\la^{(s)}}(b_s)$ and an arbitrary
differential operator $\D\in\Dl_\Llb$.
\end{rem}

\section{Construction of an eigenvector from a differential operator}

\subsection{Epimorphism $F_{\bs\la}$}
\label{Construction of B hom}
Let $\bs\la^{(1)},\dots,\bs\la^{(k)},\>\bs\la$ \,be partitions with at most $N$
parts such that $|\bs\la|=\sum_{s=1}^k|\bs\la^{(s)}|$, and $b_1,\dots,b_k$
distinct complex numbers. Denote \,$n=|\bs\la|$ \,and \,$n_s=\>|\bs\la^{(s)}|$,
\,$s=1,\dots,k$.

For $s=1,\dots,k$, \,let \,$F_s:V^{\otimes n_s}\to L_{\bs\la^{(s)}}$
be an epimorphism of $\gln$-modules. Then
\beq
\label{FF}
F_1\otimes\dots\otimes F_k:\otimes_{s=1}^kV(b_s)^{\otimes n_s}\,\to\,
\otimes_{s=1}^k L_{\bs\la^{(s)}}(b_s)
\eeq
is an epimorphism of $\glNt$-module, which induces an epimorphism
of $\B$-modules
\bea
\label{B homo}
F:(\otimes_{s=1}^kV(b_s)^{\otimes n_s})_{\bs\la}^\sing\,\to\,
(\otimes_{s=1}^k L_{\bs\la^{(s)}}(b_s))_{\bs\la}^\sing\,.
\eea

\subsection{Main result}
\label{Main result}
Let $\D^0$ be an element of $\Dl_\Llb$. Let $X^0$ be the kernel of $\D^0$.
Then $X^0$ is a point of the cell $\Oll$. Choose a germ of an algebraic curve
$X(\ep)$ in $\Oll$ such that $X(0)=X^0$ and $X(\ep)$ are generic points of
$\Oll$ for all nonzero $\ep$. Let $\bs t(\ep)$ be the root coordinates of
$X(\ep)$. The algebraic functions $t_1^{(0)}(\ep),\dots,t_n^{(0)}(\ep)$ are
determined up to permutation. Order them in such a way that the first $n_1$ of
them tend to $b_1$ as $\ep\to 0$, the next $n_2$ coordinates tend to $b_2$, and
so on until the last $n_k$ coordinates tend to $b_k$.

For every nonzero $\ep$, the vector $v(\ep)=\omega(\bs t(\ep))$ belongs to
$(V^{\otimes n})_{\bs\la}^\sing$. This vector is an eigenvector of the Bethe
algebra $\B$, acting on $(\otimes_{s=1}^nV(t_s^{(0)}(\ep)))_{\bs\la}^\sing$,
and we have $\D^{\mc B}_{v(\ep)}=\D_{X(\ep)}$, see Theorem~\ref{thm X to Vn}.

The vector $v(\ep)$ algebraically depends on $\ep$.
Let $v(\ep)=v_0\>\ep^{a_0}+v_1\>\ep^{a_1}+\dots{}$ be its Puiseux expansion,
where $v_0$ is the leading nonzero coefficient.

\begin{thm}
\label{main thm}
For a generic choice of the maps $F_1,\dots,F_k$, the vector\/ $F(v_0)$ is
nonzero. Moreover, $F(v_0)$ is an eigenvector of the Bethe algebra\/ \>$\B$,
acting on $(\otimes_{s=1}^kL_{\bs\la^{(s)}}(b_s))_{\bs\la}^\sing$, and\/
$\D^\B_{F(v_0)}=\D^0$.
\end{thm}

\begin{proof}
For any element $B\in\B$, the action of $B$ on the $\Uglnt$-module
$\otimes_{s=1}^nV(z_s)$ determines an element of $\End(V^{\otimes n})$,
polynomially depending on $z_1,\dots,z_n$. Since for every nonzero $\ep$,
the vector $v(\ep)$ is an eigenvector of $\B$, acting on
$(\otimes_{s=1}^nV(t_s^{(0)}(\ep)))_{\bs\la}^\sing$, and
$\D^{\mc B}_{v(\ep)}=\D_{X(\ep)}$, the vector $v_0$ is an eigenvector of $\B$,
acting on $(\otimes_{s=1}^kV(b_s)^{\otimes n_s})_{\bs\la}^\sing$, and
$\D^{\mc B}_{v_0}=\D^0$.

The $\glnt$-module $\otimes_{s=1}^kV(b_s)^{\otimes n_s}$ is a direct sum of
irreducible $\glnt$-modules of the form $\otimes_{s=1}^k L_{\bs\mu^{(s)}}$,
where $|\bs\mu^{(s)}|=n_s$, \,$s=1,\dots,k$. Since $\D^0\in\Dl_\Llb$,
the vector $v_0$ belongs to the component of type
$\otimes_{s=1}^kL_{\bs\la^{(s)}}(b_s)$. Therefore, for a generic choice
of the maps $F_1,\dots,F_k$, the vector\/ $F(v_0)$ is nonzero.

Since the map $F_1\otimes\dots\otimes F_k$, see~\Ref{FF},
is a homomorphism of $\glnt$-modules, the vector $F(v_0)$ is
an eigenvector of the Bethe algebra \>$\B$, acting on
$(\otimes_{s=1}^kL_{\bs\la^{(s)}}(b_s))_{\bs\la}^\sing$, and
\>$\D^\B_{F(v_0)}=\D^0$.
\end{proof}

\begin{rem}
The direction of the vector $v_0$ can depend on the choice of the algebraic
curve $X(\ep)$ in $\Om_{\bs\la}$. However, Theorem~\ref{thm eigen to Fuchs}
yields that the direction of the vector $F(v_0)$ does not depend on either
the choice of the curve $X(\ep)$ or the choice of the maps $F_1,\dots,F_k$.
\end{rem}

Given $\D\in\Dl_\Llb$\>, denote by \,$w(\D)$ the vector
\>$F(v_0)\in(\otimes_{s=1}^kL_{\bs\la^{(s)}}(b_s))_{\bs\la}^\sing$ constructed
from $\D$ in Section~\ref{Main result}. The vector $w(\D)$ is defined up to
multiplication by a nonzero number. The assignment \>$\D\mapsto w(\D)$ gives
the correspondence, which is inverse to the correspondence \>$v\mapsto\D^\B_v$
\,in Theorem~\ref{thm eigen to Fuchs}.

\subsection{Completeness of Bethe ansatz for $\glN$ Gaudin model}
\label{Sec Bethe ansatz}

The construction of the vector
$w(\D)\in(\otimes_{s=1}^kL_{\bs\la^{(s)}}(b_s))_{\bs\la}^\sing$ from
a differential operator $\D\in\Dl_\Llb$ can be viewed as a (generalized)
Bethe ansatz construction for the $\glN$ Gaudin model, cf.~the Bethe ansatz
constructions in~\cite{Ba}, \cite{RV}, \cite{MV1}, \cite{MV2}.

The following statement is contained in Theorem~6.1, Corollary~6.2 and
Corollary~6.3 of~\cite{MTV4}.

\begin{thm}
\label{thm simplicity}
If\/ $b_1,\dots,b_k$ are distinct real numbers, then the action of the Bethe
algebra on $(\otimes_{s=1}^kL_{\bs\la^{(s)}}(b_s))_{\bs\la}^\sing$ is
diagonalizable and has simple spectrum.
\end{thm}

Hence, for generic complex numbers $b_1,\dots,b_k$, there exists an eigenbasis
of the action of the Bethe algebra on
$(\otimes_{s=1}^kL_{\bs\la^{(s)}}(b_s))_{\bs\la}^\sing$. This eigenbasis is
unique up to permutation of vectors and multiplication of vectors by nonzero
numbers.

\begin{cor}
\label{cor 1}
If\/ $b_1,\dots,b_k$ are distinct real numbers or $b_1,\dots,b_k$ are generic
complex numbers, then the collection of vectors
\vvn.4>
\be
\{\>w(\D)\in (\otimes_{s=1}^kL_{\bs\la^{(s)}}(b_s))_{\bs\la}^\sing
\ |\ \D\in\Dl_\Llb\>\}
\vv.4>
\ee
is an eigenbasis of the action of the Bethe algebra.
\end{cor}

\end{document}